\documentclass[11pt]{amsart}

\usepackage{amsmath}
\usepackage{amssymb}
\usepackage{enumerate}

\begin{document}

\makeatletter
\makeatother

\numberwithin{equation}{section}

\theoremstyle{plain}
\newtheorem{thm}{Theorem}[section]
\newtheorem{prop}[thm]{Proposition}
\newtheorem{lemma}[thm]{Lemma}
\newtheorem{cor}[thm]{Corollary}
\theoremstyle{definition}
\newtheorem{remark}[thm]{Remark}
\newtheorem{remarks}[thm]{Remarks}

\newtheorem{THEO}{Theorem\!\!}
\newtheorem{DEFI}{Definition\!\!}
\renewcommand{\theTHEO}{}
\renewcommand{\theDEFI}{}

\newenvironment{ass}[1]
{\vskip 6pt \noindent {\textbf{H#1}} \rightskip 1pt \it }
{\par\vskip 6pt}

\newcommand{\refthm}[1]{Theorem~\ref{#1}}
\newcommand{\refprop}[1]{Proposition~\ref{#1}}
\newcommand{\reflemma}[1]{Lemma~\ref{#1}}
\newcommand{\refremark}[1]{Remark~\ref{#1}}

\newcommand{\fin}{\vspace{-0.3cm}
                  \begin{flushright}
                  \mbox{$\Box$}
                  \end{flushright}
                  \noindent}

\renewcommand{\proofname}{Proof}

\def\aa{\alpha}
\def\bL{{\bf L}}
\def\CC{{\mathcal C}}
\def\cL{{\mathcal L}}
\def\cR{{\mathcal R}}
\def\cVa{{\mathcal V}_\aa}
\def\d{\, \mathrm{d}}
\def\Da{D_\aa}
\def\DD{{\mathcal D}}
\def\DDDm{{\tilde\DD}_+^\aa}
\def\DDDp{{\tilde\DD}_-^\aa}
\def\DDa{\DD_\aa}
\def\DDm{\DD_+^\aa}
\def\DDp{\DD_-^\aa}
\def\Dm{{\rm Dom}}
\def\Dcma{{}^C \! D^\aa_{+}}
\def\Dpa{D^\aa_{-}}
\def\Dma{D^\aa_{+}}
\def\Dmma{\Delta^\aa_+}
\def\Dmam{D^{\aa-1}_{+}}
\def\Ea{E_\aa}
\def\Eaa{{\tilde E_\aa^q}}
\def\Fa{F_\aa}
\def\EE{{\mathbb E}}
\def\hEE{{\hat \EE}}
\def\elaw{\stackrel{{\rm d}}{=}}
\def\Jam{J_-}
\def\Jap{J_+}
\def\hbL{{\bf {\hat L}}}
\def\hpsi{{\hat \psi}}
\def\hL{{\hat L}}
\def\hR{{\hat R}}
\def\hS{{\hat S}}
\def\hT{{\hat T}}
\def\hU{{\hat U}}
\def\huu{{\hat u}}
\def\hX{{\hat X}}
\def\hxi{{\hat \xi}}
\def\hY{{\hat Y}}
\def\htY{}
\def\hZ{{\hat Z}}
\def\tZ{{\tilde Z}}
\def\tY{{\tilde Y}}
\def\ii{{\rm i}}
\def\lbd{{\lambda}}
\def\Pam{P^{-}}
\def\Pap{P^{+}}
\def\pb{{\mathbb P}}
\def\hpb{{\hat \pb}}
\def\QQ{{\mathbb Q}}
\def\hcR{{\hat \cR}}
\def\psim{{\psi_{+}}}
\def\psip{{\psi_{-}}}
\def\r{{\mathbb R}}
\def\Ta{{T_\aa}}
\def\thm{{\theta_+}}
\def\thp{{\theta_-}}
\def\Tm{{T_{-}}}
\def\Tp{{T_{+}}}
\def\Uam{U^{-}}
\def\Uamm{{\bf {\hat U}}}
\def\Uap{U^{+}}
\def\Uapp{{\bf U}}
\def\Un{{\bf 1}}
\def\Va{V_\aa}
\def\xim{{\xi^{+}}}

\def\lacc{\left\{}
\def\lcr{\left[}
\def\lpa{\left(}
\def\lva{\left|}
\def\racc{\right\}}
\def\rcr{\right]}
\def\rpa{\right)}
\def\rva{\right|}

\title[Intertwining fractional derivatives]{Intertwining certain fractional derivatives}

\author[P. Patie]{Pierre Patie}
\address{D\'epartement de Math\'ematiques, Universit\'e Libre de Bruxelles\\
Boulevard du Triomphe,  B-1050 Bruxelles. {\em Email}: {\tt ppatie@ac.ulb.be}}
\author[T. Simon]{Thomas Simon}
 \address{Laboratoire Paul Painlev\'e, U. F. R. de Math\'ematiques, Universit\'e de Lille 1, F-59655 Villeneuve d'Ascq Cedex. {\em Email}: {\tt simon@math.univ-lille1.fr}}



\bibliographystyle{plain}
\begin{abstract}
We obtain an intertwining relation between some Riemann-Liouville operators of
order $\aa \in (1,2),$ connecting through a certain multiplicative identity in
law the one-dimensional marginals of reflected completely asymmetric
$\aa-$stable L\'evy processes. An alternative approach based on recurrent extensions of positive self-similar Markov processes and exponential functionals of L\'evy processes is also discussed.
\end{abstract}

\maketitle

\section{Introduction and statement of the result}

Consider for every $\aa \in (1,2)$ the following operators acting on functions from $\r^+$ to $\r$:
$$\Dmam f(x) \;=\; \frac{\d}{\d x}\int_0^{x}
\frac{f(t)(x-t)^{1-\alpha}}{\Gamma(2-\alpha)}\, dt\qquad\mbox{and}\qquad
\Dpa f(x) \;=\; \frac{\d^2}{\d x^2}\int_x^{\infty}
\frac{f(t)(t-x)^{1-\alpha}}{\Gamma(2-\alpha)}\, dt,$$
for every $x >0.$ The operator $\Dmam$ is known as the left-sided   Riemann-Liouville (RL)
derivative of index $\aa-1$ and $\Dpa$ as the right-sided RL derivative of
index $\aa$. Recently, RL derivatives have appeared quite often in various
domains of analysis and probability and we refer to Chapter 2 in \cite{KST}
for a detailed account on them, as well as on other fractional
operators. From the analytical viewpoint, RL derivatives extend in a non-local fashion the
derivatives of integer order - see (2.2.5) in \cite{KST}. To name but one classical example of the
occurence of RL derivatives in probability, recall that
$$\Dpa\; =\; D^{2}\;\circ\; U^{2-\aa}$$
where $D$ is the usual derivative and $U^{2-\aa}$ the potential associated
to the standard $(2-\aa)$-stable subordinator - see e.g. Exercise (1.6) in \cite{BG}. Consider now the operator
$$\Dmma \; =\; \Dmam\,\circ\; D,$$
noticing first that it differs from $D\, \circ\, \Dmam \,=\, \Dma$ - see
(2.2.29) in \cite{KST} for this latter equality - where $\Dma$ is the
left-sided RL derivative of index $\aa$ which is analogously defined by
$$\Dma f(x) \;=\; \frac{\d^2}{\d x^2}\int_0^x \frac{f(t)(x-t)^{1-\alpha}}{\Gamma(2-\alpha)}dt.$$
Indeed, an integration by parts shows that for every $x > 0$
\begin{eqnarray}
\label{CRL}
\Dma f(x) & = & \Dmma f(x) \;+\; \frac{x^{-\alpha}}{ \Gamma(1-\alpha)}f(0).
\end{eqnarray}
Another fractional operator related to $\Dmma$ is the so-called Caputo
$\aa$-fractional derivative which is given by
$$\Dcma f(x)\; =\;\frac{1}{\Gamma(2-\alpha)}\int_{0}^xf''(u)(x-u)^{1-\alpha}du,$$
see section 2.4 in \cite{KST}. A similar integration by parts shows namely that
\begin{eqnarray}
\label{Cap}
\Dmma f(x) & = & \Dcma f(x)\;+\; \frac{x^{1-\alpha}}{ \Gamma(2-\alpha)}f'(0).
\end{eqnarray}
For any $\aa\in(1,2),$ let $\cVa$ be the multiplicative kernel acting on
functions from $\r^+$ to $\r$ by
$$\cVa f(x) \;=\; \EE\lcr f(x\Va)\rcr,$$
where $\Va$ is a positive random variable having the density
$$v_{\alpha}(t)\; =\;\frac{ (-\sin\pi\aa)t^{\alpha-2}(1+t)}{\pi(t^{2\alpha}-2t^{\alpha}\cos \alpha \pi+1)}$$
- it will be checked soon afterwards that $v_\aa$ is indeed a density
function on $\r^+$. Setting $\CC^2_b$ for the set of twice continuously
differentiable functions $\r^+\to\r$ such that $f'$ and $f''$ are bounded, consider finally the following domain
$$\DD\; =\; \{ f\in\CC^2_b\;\;\mbox{such that}\;\; f'(0) = 0\;\mbox{and}\;
\exists \; \gamma > 2-\aa \;\; /\; \lim_{x\to+\infty} x^\gamma (\vert f(x)\vert
+\vert f''(x))\vert) = 0\}.$$
Observe that if $f\in\DD$, then necessarily $\exists \; \gamma > 2-\aa \; /\;\lim_{x\to+\infty} x^\gamma \vert f'(x)\vert
 = 0$ as well, which can be checked from the decomposition
$$f(x+1) \; =\; f(x)\; +\; f'(x)\; +\; \int_x^{x+1} (x+1-y) f''(y) dy.$$
It is also easily seen that $\Dpa$ and $\cVa$ are well-defined on $\DD,$ and it will
be proved in the next section that $\Dmma$ is well-defined on $\cVa(\DD).$ Our
main result is the following intertwining relation
between $\Dmma$ and $\Dpa$:

\begin{THEO} For any $f\in\DD,$ one has
\begin{eqnarray*}
\Dmma \cVa f &=& \cVa\Dpa f.
\end{eqnarray*}
\end{THEO}

Intertwining relations between Markov processes have some history and we refer
to \cite{CPY2, HY} for a probabilistic account,  as well as various examples and
applications. See also \cite{M} for a particular analytical study. Though
expressed in analytical terms, our result falls within the Markovian
framework. It was namely observed in \cite{BDP} that $\Dmma$ (resp.~$\Dpa$) is
the infinitesimal generator of the spectrally positive (resp. spectrally
negative) $\aa-$stable L\'evy process reflected at its running  supremum. It
might be interesting to mention that although the underlying stable L\'evy
processes are in classical duality, this is no more the case for the reflected
processes, so that our result can be viewed as a kind of intertwisted duality
relationship for the latter.

The proof of the theorem hinges upon the well-known criterion given in the Proposition 3.2 of \cite{CPY1} and an identity in law connecting the running suprema of completely asymmetric stable processes which was recently obtained in \cite{ThS1}, involving the positive variable $\Ta$ with  density $$\frac{(-\sin\pi\aa)(1+ t^{1/\aa})}{\pi\aa(t^2 -2t \cos \alpha \pi+1)}\cdot$$
Notice in passing that our above function $v_\aa$ is the density of the
variable $T_\aa^{-1/\aa},$ hence it is a density function. Some particular
attention is paid to the functional domain upon which the intertwining
relation holds and the set $\DD,$ which we borrowed from \cite{BDP}, appears
to be a reasonable and not too small candidate. With the help of some
Suprun-type formul\ae \, for the resolvent of spectrally one-sided L\'evy
processes which had been derived in \cite{Pi}, we can also identify the
Fellerian domains of  $\Dmma$ and $\Dpa.$ In theory, those Fellerian domains
yield an optimal formulation for the theorem, although they do not seem
very tractable.

In Section 3 we discuss another approach which consists in interpreting the
stable process reflected at its infimum as the unique self-similar recurrent extension
leaving 0 continuously of the stable process killed when it enters the
negative half-line. This identification had been roughly explained in Example
3 of \cite{R} and here we can also check analytically that these two Feller
processes have the same infinitesimal generators. This gives another proof,
looking somewhat more unified, of all the results contained in the Appendix in
\cite{BDP}. An independent proof
of the identity in law between suprema of completely asymmetric stable
processes, which is the key-argument for the theorem, is also proposed, involving some closed formul\ae\, for the exponential functionals of certain L\'evy processes which had been established in \cite{Pa3, Pa4}. Though overall a bit lenghtier, we believe that this second point of view provides some unity to our interweaving relationship, which appears to be coherent with several apparently disconnected identities.

To conclude this introduction, we stress that the positive random variable $Z_\beta$ with density
$$\frac{(-\sin\pi\aa)}{\pi(\aa-1)(t^2 -2t \cos \alpha \pi+1)}\;=\; \frac{\sin\pi\beta}{\pi\beta (t^2 +2t \cos \beta \pi+1)}$$
where $\beta = \aa -1\in(0,1),$ which can be viewed as a cut off Cauchy
variable, has already occured in several distinct areas of the literature,
especially through its power transforms $Y_\beta = Z_\beta^{1/\beta}.$ See for instance Formula XI.11.6 in \cite{Y} for connections with the
$\beta$-fractional power of linear operators,  Theorem 1.1 and Theorem 1.2 in
\cite{Ko} for mixture representations of the Linnik and Mittag-Leffler distributions of index $\beta$, and more general geometric stable
distributions, or Exercise 4.21 (3) in \cite{CY} which shows that $Y_\beta$ has the same law as the independent quotient of two standard positive $\beta-$stable laws. The above variables $\Va$ and $\Ta$ are less classical than $Z_\beta$ but one may of course wonder if they are not particular instances of a family of positive variables connecting suprema of general stable processes in duality, or a broader class of fractional operators than the one we consider in the present article.  We plan to tackle this question in some further research.

\section{Proof of the theorem}

Let $(Z, \pb)$ be a spectrally negative stable L\'evy process of index $\alpha \in
(1,2)$ starting at $0,$ with L\'evy density
$$\nu_\aa(y) \;=\;  \frac{|y|^{-(\alpha+1)}}{\Gamma(-\alpha)}\Un_{\{ y<0\}},$$
so that the L\'evy-Khintchine formula reads
\begin{eqnarray}
\label{norm}
\EE\lcr e^{\lbd Z_t}\rcr \; =\; e^{t\lbd^\aa}
\end{eqnarray}
for every $t, \lbd \ge 0.$ See e.g. Chapters VII \& VIII in \cite{B} for an account
on completely asymmetric stable processes. Setting $S_t = \sup\{ Z_s, \,
s\le t\}$ for the associated running supremum and introducing
$$T_x =\inf\{ t > 0, Z_t = x\}
=\inf\{ t > 0, S_t = x\},$$
recall that $T_1 \elaw S_1^{-\aa}$ is a standard
positive $(1/\aa)$-stable law, viz.
$$\EE\lcr e^{\lbd T_1}\rcr \; =\; e^{-\lbd^{1/\aa}}$$
for every $\lbd \ge 0.$ Denoting by $I_t = \inf\{Z_s, \; s\le t\}$ the running
infimum, consider the reflected processes
$$X_t\; =\; S_t\; -\; Z_t\quad\mbox{and}\quad \hX_t\; =\; Z_t\; -\; I_t.$$
Notice that if $\hZ = -Z$ is the dual process and if $Y_t = \hS_t - \hZ_t$ and $\hY_t =
\hZ_t - {\hat I}_t$ are the corresponding reflected processes, then $Y = \hX$ and $\hY =
X.$ It is a basic fact from fluctuation theory - see e.g. Proposition VI.1 in
\cite{B} - that $X$ and $\hX$ are
Feller processes and we will denote by $\pb_x$ resp. $\hpb_x$ their laws
starting from $x\ge 0.$ The infinitesimal generators $L$ and $\hL$ of $X$ and $\hX$ have
recently been expressed in \cite{BDP}, in three different forms. Recalling that by definition
\begin{equation}
\label{GI}
L f(x)\; =\; \lim_{t\to 0}\frac{\EE_x[f(X_t)] -
  f(x)}{t}\qquad\mbox{resp.}\qquad  \hL f(x)\; =\; \lim_{t\to 0}\frac{\hEE_x[f(\hX_t)] -
  f(x)}{t}
\end{equation}
for every continuous function $f : \r^+ \to \r$ such that the limit in the right-hand side of the first
(resp. second) equality exists uniformly, let us denote by Dom $L$ (resp.~Dom $\hL$) the
set of such functions. Choosing then Riemann Liouville's form in Proposition A.1 of
\cite{BDP}, one has $\DD  \subset$ Dom $L$ and  
$$L f(x)\; =\; \Dpa f(x), \quad x > 0$$
for every $f\in \DD$. Besides, one has 
$$\hL f (x) \; =\; \Dmma f (x), \quad x > 0 $$
for every $f\in \CC^2_b$ such that $f'(0) = 0.$ The next proposition shows that $\Dmma$ is well-defined on
$\cVa(\DD)$ so that the statement of our theorem makes sense:

\begin{prop}
\label{SMS}
One has $\cVa (\DD)\subset$ {\rm Dom} $\hL.$
\end{prop}

\begin{proof} Since $v_\aa (t)$ is of order $t^{-(1+\aa)},$ one sees by
  dominated convergence that for any $f\in\DD$, the function $\cVa f$ is
  continuously differentiable on $\r^+$ with bounded derivative
$$(\cVa f)'(x)\; =\; \EE[ \Va f'(x\Va)]$$
- whence in particular $(\cVa f)'(0) = 0$, and twice continuously
differentiable on $(0, +\infty)$ with second derivative
$$(\cVa f)''(x)\; =\; \EE[ \Va^2 f''(x\Va)], \quad x >0.$$
On the other hand, the right-hand side in the above equality might not be bounded when $x\to 0$ because $\EE [\Va^2] = +\infty.$ An easy change of variable shows however that for every $f\in\DD,$ the quantity
$$(x^\gamma\wedge 1)\int_0^\infty \!(t^{1-\aa}\wedge 1) \, f''(xt) dt$$
remains bounded on $(0, +\infty)$ for some $\gamma < 2-\aa,$ so that $(x^\gamma\wedge 1)
(\cVa f)''(x)$ is bounded on $(0, +\infty)$. Hence, we need to show that
in Proposition 4 of \cite{BDP}, the global boundedness assumption on
$F''$ can be relaxed into $(x^\gamma\wedge 1)F''(x)$ bounded on $(0, +\infty)$ for some $\gamma\in [0, 2-\aa).$ A persual of the proof shows that it is indeed the case: the relaxed
condition changes namely the right-hand side of the second inequality in
(4.11) therein into
$$C_\aa{\bf E}_y\lcr \int_0^t (Y_s^{2 -\aa} + Y_s^{2-\aa-\gamma})\, ds \rcr$$
for some positive finite constant $C_\aa,$ where ${\bf P}_y$ stands for the law of $Y$ (with their notation for $Y,$ which matches ours) starting at $y \ge 0.$ Using self-similarity, the expectation is then bounded by
$$C'_\aa ( {\bf E}_y[Y_1^{2-\aa}] + {\bf E}_y[Y_1^{2 -\aa-\gamma}])\; \le\; C_{\aa, \gamma} {\bf E}_y[Y_1^{2-\aa}]$$
for some other positive finite constants $C_\aa', C_{\aa, \gamma}$. Notice now by the Feller property and the identification $Y = \hX$ that the law of $Y$ under ${\bf P}_y$ is that of the process $\{Z_t + (y\vee (-I_t)), \; t\ge 0\}$ under $\pb.$ Since $0 <
2-\aa < 1 <\aa,$ one has e.g. from Proposition VIII.4 in \cite{B} and the
comment thereafter
$${\bf E}_y[Y_1^{2-\aa}]\; =\; \EE[S_1^{2-\aa}]\; +\;  \EE[(y \vee
(-I_1))^{2-\aa}]\; <\; +\infty,$$
which was already used in \cite{BDP}. All of this shows that (4.12) remains unchanged for $F$ under the relaxed condition, which entails exactly as in \cite{BDP} that the limit in the right-hand side of (\ref{GI}) exists pointwise for any $f\in\cVa(\DD).$  To finish the proof, notice first by the Feller property that the transition densities of $\hX$ form a strongly continuous contraction semigroup on the Banach space of continuous functions $\r^+\to\r$ tending to 0 at infinity, to which belongs $\cVa(\DD).$ Hence, one can apply the whole semigroup theory recalled e.g. in Section 31 of \cite{S}. In particular, a result of K.~It\^o - see Lemma 31.7 in \cite{S} shows that the limit in the right-hand side of (\ref{GI}) is actually uniform for any $f\in\cVa(\DD),$ whence $f\in{\rm Dom}\;\hL$ as desired.

\end{proof}

Let us now consider the Fellerian domains of $L$ and $\hL.$ Setting $\CC_0$ for the set of continuous functions $\r^+\to\r$ tending to 0 at infinity, the Feller property states that both functions
$$x\;\mapsto\; \EE_x[f(X_t)]\quad\mbox{and}\quad x\;\mapsto\;
\hEE_x[f(\hX_t)]$$
are in $\CC_0$ whenever $f\in\CC_0.$ The Fellerian domains of $L$ and $\hL,$ which we denote respectively by $\DD(L)$ and $\DD(\hL)$ are made out of those functions in
$\CC_0$ such that the limits in (\ref{GI}) exist uniformly - see
e.g. Definition VII (1.1) in \cite{RY}. It is
well-known from semigroup theory - see e.g. Proposition
VII (1.4) in \cite{RY} - that
\begin{equation}
\label{SG}
\DD (L) \; =\;  U_q(\CC_0)\;\;\mbox{and}\;\; \DD (\hL) \;=\; \hU_q(\CC_0)
\end{equation}
for every $q > 0,$ where $U_q, \hU_q$ are the resolvent operators defined by
$$U_q f(x) \; =\; \int_0^\infty e^{-qt} \EE_x[f(X_t)]\, dt\quad \mbox{and}
\quad \hU_q f(x) \; =\; \int_0^\infty e^{-qt} \hEE_x[f(\hX_t)]\, dt.$$
It follows from the definitions that
$${\rm Dom }\, L \;\cap\;\CC_0\; =\; \DD(L)\quad\mbox{and}\quad {\rm Dom }\, \hL
\;\cap\;\CC_0\; =\; \DD(\hL).$$
The next proposition gives a full description of $\DD(L)$ and $\DD(\hL)$ in
terms of the functions $\Fa (x) = \Ea(x^\aa)$ and its derivatives, where $\Ea$ is the
Mittag-Leffler function of index $\aa$ which is defined by
$$\Ea(x)\;=\;\sum_{n=0}^{\infty} \frac{x^n}{\Gamma(\alpha n+1)}\cdot$$
For every $f\in\CC_0,$ introduce the further notation
$$\lbd_f\; =\;\int_0^\infty e^{-y} f(y) dy.$$
\begin{prop}
\label{Dom}
One has $\DD(\hL)\; =\; \lacc x\mapsto  \lbd_f \Fa(x)\, -\,\Fa'\star f \, (x), \; \; f\in\CC_0\racc$
and
$$\DD(L)\; =\; \lacc x\mapsto  e^{-x}\int_0^x \Fa''(y) f(y) dy \,
+\,\int_x^\infty (e^{-x}\Fa''(y) - \Fa'(y-x)) f(y) dy, \; \; f\in\CC_0\racc.$$
\end{prop}

\begin{proof} Let us start with $\DD(\hL) = \hU^1(\CC_0).$ The resolvent density of $\hX = Y$ killed when entering the half-line  $(a,\infty)$ has been computed for every $a > 0$ in \cite{Pi}, Theorem 1 (i). Notice first that the notations therein yield $Z^{(1)}(x) =\Fa (x)$ and $W^{(1)}(x) = \Fa'(x).$ Besides, it was proved in
  Theorem 1 of \cite{ThS1} that the function $\Fa - \Fa'$ is completely monotone, so that
  in particular $\Fa (x) -\Fa'(x) \to 0$ as $x\to +\infty$ which, together with
  Formula 18.1 (10) in \cite{E}, entails
$$\frac{W^{(1)}(a-y)}{Z^{(1)}(a)}\; \to \; e^{-y}\quad\mbox{as $a\to +\infty$}$$
for every $y\ge 0.$ By monotone convergence, letting $a \to +\infty$ shows that the density of $\hU^1$ is
$$\huu^1(x,y)\; =\; e^{-y} \Fa(x) \; -\; \Fa'(x-y)\Un_{\{y\le x\}},$$
whence
$$\hU^1\! f(x) \; = \; \lpa\int_0^\infty e^{-y} f(y) dy\rpa \Fa(x)\; -\; \int_0^x
\Fa'(x-y) f(y) dy \; =\;  \lbd_f \Fa(x)\, -\,\Fa'\star f \, (x)$$
for every $f\in\CC_0$ and $x\ge 0$, which is the required identification for $\DD(\hL).$

Before identifying $\DD(L) =\hU^1(\CC_0),$ let us check that the two integrals in the definition are well-defined for every $f\in\CC_0$ and $x\ge 0:$ the convergence of the proper integral comes from the easily shown behaviour $\Fa''(y)
\sim (\aa- 1) y^{\aa-2}/\Gamma(2\aa)$ as $y\to 0,$ and the existence of the
improper one is proved in combining several times formul\ae\, (6) and (43) in
\cite{GLL}, which yield
$$\aa\Fa'(z) \; =\; e^z\; +\; {\rm O} (
z^{-(1+\aa)})\quad\mbox{and}\quad \aa\Fa''(z) \; =\; e^z\; +\; {\rm O} (
z^{-(2+\aa)})$$
as $z\to +\infty.$ The latter asymptotics also entail, with the notations of \cite{Pi},
$$\frac{W^{(1)}(a-x)}{W^{(1)'}_+(a)}\; =\; \frac{\Fa'(a-x)}{\Fa''(a)}\; \to \; e^{-x}\quad\mbox{as}\; a\to +\infty.$$
Since moreover $W^{(1)}(0) = \Fa'(0) = 0,$ one
obtains from Theorem 1 (ii) in \cite{Pi} and monotone convergence the following
expression for the density of $U_1:$
$$u^1(x,y)\; =\; e^{-x} \Fa''(y) \; -\; \Fa'(y-x)\Un_{\{y\ge x\}}$$
whence, as above, the desired expression for $\DD(\hL).$

\end{proof}

\begin{remarks} (a) If $g = \hU^1\! f \in\DD (\hL),$ then $g(0) = \lbd_f$ and $g$ is
  continuously differentiable on $\r^+$ with derivative $g'(x) = \lbd_f
  \Fa'(x) - \Fa''\star f(x),$ so that in particular $g'(0) = 0.$ On the other
  hand, choosing $f$ non differentiable shows that $\DD (L)$ contains
  functions which are not $\CC^2,$ as of course might be expected from the
  expression of $L.$ See Chapter 2 in \cite{KST} for more material concerning
  the domains of Riemann-Liouville derivatives.

\vspace{2mm}

(b) If $g = U^1\! f \in\DD (L),$ then
$$g(0)\; = \; \int_0^\infty (\Fa''(y) -\Fa'(y)) f(y) dy$$
and $g$ is continuously differentiable on $\r^+$ with derivative
$$g'(x)\; = \; \int_x^\infty (\Fa''(y-x) - e^{-x}\Fa''(y)) f(y) dy \; -\;
e^{-x} \int_0^x \Fa''(y) dy,$$
so that in particular $g'(0) = 0.$ Again,  $\DD (\hL)$ contains
  functions which are not $\CC^2.$

\end{remarks}

\noindent
{\bf End of the proof.} We will use the inclusion (2) $\Rightarrow$ (3) in the
Proposition 3.2 of \cite{CPY1}. We already know that $X$ and $\hX$ are Feller
processes and it follows easily from the $(1/\aa)-$self-similarity of $Z$ that
they are also $(1/\aa)-$self-similar Markov processes. More precisely, one has
$$(X^b, \pb_x)\; \elaw\; (X, \pb_{b x})\quad\mbox{and}\quad (\hX^b, \hpb_x)\; \elaw\; (\hX, \hpb_{b x})$$
for every $b >0,$ with $X^b_t = b X_{b^{-\aa} t}$ and an analogous notation
for $\hX^b.$ Since $\aa\neq 1,$ this shows that they are not semi-stable
Markov processes viz.~1-self-similar Markov processes - see (1.b) in
\cite{CPY1}. However, a perusal of the proof of Proposition 3.2 in \cite{CPY1}
shows that its statement remains unchanged when considering
$(1/\aa)-$self-similar Markov processes for every $\aa > 0$ and not just $\aa =1.$

We next show that the distribution of $\hX$ is determinate under $\hpb_0,$
with the notation of \cite{CPY1}. From e.g. Proposition VI.3 in \cite{B}, the
law of $\hX_t$ under $\hpb_0$ is that of $S_t$ under $\pb.$ Hence, taking the
Laplace transforms, we need to show that if $f, g\in\CC_0$ are such that
\begin{equation}
\label{Deter}
\EE\lcr \int_0^\infty e^{-\lbd t} f(S_t) \,dt \rcr\; =\; \EE\lcr \int_0^\infty
e^{-\lbd t} g(S_t) \, dt \rcr
\end{equation}
for every $\lbd > 0,$ then $f = g.$ The latter is a basic property of $(1/\aa)-$stable
subordinators but we will give some details for the reader's
comfort. Recalling the notation $T_x = \inf\{t > 0, \; Z_t = x\}$ for every
$x\ge 0,$ we know that  $\lacc T_x, \; x\ge 0\racc$ is a
standard $(1/\aa)$-stable subordinator since $Z$ has no positive
jumps. Besides, one has ${\rm Leb}\, (\r^+ -\cup_{x\ge 0} (T_{x-}, T_x)) = 0$
a.s. - see e.g. the
beginning of Section III.5 in \cite{B}, so that (\ref{Deter}) entails
$$\EE\lcr\sum_{x>0} f(x) e^{-\lbd T_{x-}}(1 - e^{-\lbd \Delta T_x }) \rcr\;
=\; \EE\lcr \sum_{x>0} g(x) e^{-\lbd T_{x-}}(1 - e^{-\lbd \Delta T_x })\rcr.$$
The so-called Master's Formula - see e.g. Proposition XII (1.10) in \cite{RY} - yields then
$$\EE\lcr\int_0^\infty f(x) e^{-\lbd T_{x-}}\lpa\int_0^\infty
\frac{1 - e^{-\lbd u}}{u^{\aa +1}}du\rpa dx\rcr\;
=\; \EE\lcr \int_0^\infty g(x) e^{-\lbd T_{x-}}\lpa\int_0^\infty
\frac{1 - e^{-\lbd u}}{u^{\aa +1}}du\rpa dx\rcr,$$
whence
$$\int_0^\infty f(x) e^{-\lbd^{1/\aa} x} dx\; =\; \int_0^\infty g(x)
e^{-\lbd^{1/\aa} x}dx$$
and the required identification $f=g$ by inversion of the Laplace transforms.

Last, we see from Proposition VI. 3 in \cite{B} and Formula (9) in \cite{ThS1} (beware the inverse notations) that
$$X_1\; \elaw\; \hS_1\;\elaw\; T_\aa^{-1/\aa}\, \times\, S_1\; \elaw\; \Va\,\times\, S_1\; \elaw\; \Va\,\times\,\hX_1,$$
where the identity $T_\aa^{-1/\aa}\elaw \Va$ was mentioned in the introduction and follows from a change of variable. Putting everything together, Proposition 3.2 in \cite{CPY1} shows that $\Dmma\cVa f = \cVa\Dpa f$ for every
$f\in\DD(L)$ such that $\cVa f\in \DD(\hL).$ Supposing now that $f\in \DD\subset
{\rm Dom}\, L \cap \CC_0 = \DD(L),$ Proposition \ref{SMS} shows that $\cVa f
\in {\rm Dom} \, \hL$ and it follows immediately from dominated convergence
that $\cVa f \in \CC_0$. Hence, $\cVa f \in {\rm Dom}\, \hL \cap \CC_0 =
\DD(\hL)$ and we have shown
$$\Dmma\cVa f\; = \;\cVa\Dpa f$$
for every $f\in\DD,$ as required.

\fin

\begin{remarks} (a) The interweaving relationship
$$\Dcma\cVa f\; =\; \cVa\Dpa f$$
also holds over $\DD$ because of the identity (\ref{Cap}). Actually, since
$f'(0) = 0$ for every $f\in \DD(\hL),$ the
operator $\Dcma$ coincides with the generator $\hL$ of $\hX.$

\vspace{2mm}

(b) As mentioned during the proof, Proposition 3.2 in \cite{CPY1} shows that
$$\Dmma\cVa f\; =\; \cVa\Dpa f$$
holds for every $f\in\DD(L)$ such that $\cVa f\in \DD(\hL).$ The domains $\DD(L)$ and $\DD(\hL)$ have been identified in Proposition \ref{Dom} and allow rougher functions than $\DD$ but their formulations are unfortunately not very tractable, contrary to $\DD.$

\end{remarks}

We close this section with an interesting complete monotonicity property for the Mittag-Leffler function $\Ea.$ The latter is actually a direct consequence of Proposition 2 in \cite{Pi}, but we present here a proof based on generators which is perhaps more transparent. Recall that a smooth function $f : ]0, +\infty[ \to\r^+$ is said to be completely monotone (CM) if
$$(-1)^n \frac{\d^n \! f}{\d x^n}\; \ge \; 0, \quad n\ge 1.$$
By Bernstein's theorem, when $f(x)\to 1$ as $x\to 0$ this is equivalent to the fact that $f$ is the Laplace transform of a probability measure over $\r^+.$ A classical result by Pollard states that for any $\aa\in (0,1]$ the function $\Ea(-x)$ is CM, a property which does not hold anymore if $\aa > 1$ - see \cite{Sc} for a proof of these latter facts and more on this topic. Like $1/x$, the function $x \mapsto \Ea(1/x)$ is CM by positivity of the coefficient in the series expansion - see Section 4 in \cite{MS} for further properties as well as general references on complete monotonicity. The following proposition shows a related property in the case $\aa\in [1,2]$:

\begin{prop}
\label{ML}
For any $\aa\in [1,2]$ the function $x \mapsto 1/\Ea(x)$ is {\em CM}.
\end{prop}

\begin{proof} The case $\aa = 1$ is straightforward because $E_1(x) = e^x$ and the case $\aa = 2$ follows from L\'evy's formula: one has
$$\frac{1}{E_2(x)}\; =\; \frac{1}{\cosh \sqrt{x}}\; =\; \EE[e^{-x \tau}]$$
with $\tau = \inf\{t> 0, \; \vert B_t\vert > 2^{-1/4}\}$ and $B$ a standard
Brownian motion - see e.g. Exercise II (3.10) in \cite{RY}, so that $1/E_2(x)$ is CM by Bernstein's theorem. For the
remaining case $\aa\in (1,2)$ we will follow roughly the same arguments as Theorem 2.1
and 2.6 in \cite{Pa2}. Setting
$$\Ea^q(x) \; = \;\Ea(qx^\aa)f_q(x)$$
for every $q > 0,$  where $f_q$ is some
smooth function such that $f_q(x) = 1$ for all $x\le 1$ and $f_q(x) = 0$ for
all $x\ge 2,$ say, we see that $\Ea^q\in\DD(\Dmma)$. On the other hand, it follows from Lemma 2.23 in \cite{KST} that
$$\Dmma \Ea^q(x)\; =\; q\Ea^q(x)$$
for every $x\le 1.$ From Definition VII (1.8) and Exercise VII (1.24) in \cite{RY} - this shows that $t\mapsto e^{-q T_1^+\wedge t}
\Ea^q (\hX_{T_1^+\wedge t})$ is a martingale under $\hpb_0$, where
$$T_1^+\; =\; \inf\{s>0, \; \hX_s\geq 1\}\; =\; \inf\{s>0, \; \hX_s = 1\}$$
(recall that $\hX$ has no positive jumps for the second equality). The optional sampling theorem entails
$$1\; =\; \hEE_0 [ e^{-qT_1^+\wedge t} \Ea^q (\hX_{T_1^+\!\wedge t})] \; \to \;
\hEE_0 [ e^{-qT_1^+} \Ea^q (\hX_{T_1^+})]\; =\;\Ea(q)\hEE_0 [ e^{-qT_1^+}]$$
as $t\to\infty$, by dominated convergence because $\hX_{T_1^+\!\wedge t}\in[0,1]$ a.s. This completes the proof by Bernstein's theorem.

\end{proof}

\begin{remark} It is plain from self-similarity, the Markov property and the
  absence of positive jumps that the variable $T^+_1$ is self-decomposable. By
  a well-known argument - see the proof of Proposition 2.6 in \cite{Pa2} for
  details - this shows that the mapping
$$x \;\mapsto\; \exp -[x\Ea'(x)/\Ea(x)]$$
is also CM as the Laplace transform of a positive infinitely divisible random variable.
\end{remark}

\section{Another approach with recurrent extensions}

The purpose of this section is two-fold. First, we will give another proof of the identifications $L = \Dpa$ and $\hL =\Dmma,$ viewing the reflected process $X$ (resp.~$\hX$) as a recurrent extension of the process $\hZ$ (resp.~$Z$) killed when
entering the negative half-line, and using the classical expression of the
infinitesimal generator of the unkilled stable L\'evy process. Actually the approach works for every strictly stable process $Z$ such that $\vert Z\vert$ is not a subordinator, making it possible to retrieve the whole Appendix of \cite{BDP}. Second, we will derive a proof of the identity
\begin{equation}
\label{Cruz}
X_1\; \elaw\; \Va\,\times\, \hX_1,
\end{equation}
which is independent of \cite{ThS1} and relies upon closed expressions for the densities of exponential functional of certain spectrally negative L\'evy processes that had been carried out in \cite{Pa3, Pa4}.

\subsection{Retrieving the generators of reflected stable processes} Let $(Z, \QQ_x)$ be a strictly stable L\'evy process of index $\aa\in (0,2)$ such that $\vert Z\vert$ is not a subordinator, starting from $x\in\r.$ We refer e.g. to Chapter VIII in \cite{B}  for details. The density of the L\'evy measure is
$$\nu(y) \; =\; c_+ y^{-\aa -1}\Un_{\{y > 0\}}\; +\;  c_- \vert y\vert^{-\aa -1}\Un_{\{y < 0\}},$$
where $c_+, c_-$ are nonnegative constants such that $c_+ + c_- > 0.$
When $\aa = 1$  we suppose that $Z$ is a symmetric Cauchy process viz. $c_+ =  c_- = c > 0.$ Again, we will use the notations $S_t =\sup\{Z_s, \, s\le t\}, I_t = \inf\{Z_s, \, s\le t\}, X_t = S_t - Z_t$ and $\hX_t = Z_t - I_t.$ Setting also $T =\inf\{t > 0, \; Z_t\le 0\},$
consider now the killed process
$$R_t\; =\; Z_t \Un_{\{T > t\}}.$$
Again, since $(T, \QQ_x) \elaw (x^\aa T, \QQ_1)$  for every $x\ge 0,$ one sees that $(R, \QQ_x)$ is a positive $(1/\aa)-$self-similar Markov processes, viz. a Feller process taking values in $\r^+$ (here, with 0 as an absorbing state) and fulfilling the scaling property
$$(R^b, \QQ_x)\; \elaw\; (R, \QQ_{b x})$$
for every $b >0,$ with the notation $R^b_t = b R_{b^{-\aa} t}.$
As noticed in Example 3 of \cite{R}, the reflected process $(\hX, \hpb_x)$ (with $\hpb_x$ defined
analogously as in the preceding section) can be viewed as the unique self-similar recurrent extension of $(R, \QQ_x)$ leaving 0 continuously. Roughly speaking, for every $x > 0$ those two processes have the same law until the a.s. finite time $T$ resp. $\hT = \inf\{t > 0, \; \hX_t\le 0\}$ but 0 is a regular boundary point for $\hX,$ which is left instantaneously and continuously. The Feller process $\hX$ has also infinite lifetime. See \cite{F, R} for precise accounts on recurrent extensions.

We now identify the generator $\hcR$ of the reflected process $\hX$ viewed as a recurrent extension of $R$, retrieving in a unified manner all the results contained in the Appendix of \cite{BDP}. Beware that since we consider the process reflected at its infimum, our notation is reverse to
that of \cite{BDP}. As in Section 2, we set $(\hU_q)_{q\ge 0}$ for the resolvent of $\hX.$ For every $\aa\in (0,1)$ and $f\in\DD,$ we use the same notations
$$\Dpa f(x) \; =\; \frac{\d^2}{\d x^2}\int_x^{\infty}
\frac{f(t)(t-x)^{1-\alpha}}{\Gamma(2-\alpha)}\, dt\;=\; \frac{\d}{\d x}\int_x^{\infty}
\frac{f(t)(t-x)^{-\alpha}}{\aa\Gamma(-\aa)}\, dt$$
and
$$\Dma f(x) \; =\; \frac{\d^2}{\d x^2}\int_0^x
\frac{f(t)(x-t)^{1-\alpha}}{\Gamma(2-\alpha)}\, dt\;=\; - \frac{\d}{\d x}\int_0^x
\frac{f(t)(x-t)^{-\alpha}}{\aa\Gamma(-\aa)}\, dt.$$

\begin{prop}[Bernyk-Dalang-Peskir]
\label{IDT} For every $f\in\DD$ and $x >0,$ one has
$$\hcR f(x)\; =\; \Gamma(-\aa) (c_-\Dma f (x) + c_+\Dpa f (x))\; +\; \frac{c_- f(0)}{\aa x^\aa}$$
if $\aa\neq 1,$ and
$$\hcR f (x)\; =\; c \lpa  \frac{\d^2}{\d x^2}\int_0^{\infty} f(t)\log \lpa
\frac{1}{\vert x-t\vert}\rpa dt\; +\; \frac{f(0)}{x}\rpa$$
if $\aa = 1.$
\end{prop}

\begin{proof} Fix $f\in\DD$ and $x> 0.$ For every $M > 0,$ define $f_M$ over $\r$ in setting $f_M(x) = f(x)$ for every $x\ge 0,$ $f_M(x) = f(0)$ for every $- M <x <0,$ and letting $f_M(x)\to 0$ smoothly as $x\to -\infty$. Then $f_M\in\CC^2_b(\r)$ except possibly at zero where its left and right second derivatives are bounded, and $f_M(x)\to 0$ as $\vert x\vert \to +\infty.$ Besides, with an abuse of notation, one can write $\hcR f(x) = \hcR f_M(x)$ for every $M > 0.$ Introducing the resolvent of $R$
$$\Uapp_q g(x)\; =\; \int_0^\infty e^{-qt}\; \EE_x\! \lcr g(Z_t) \Un_{\{ T >
  t\}}\rcr dt$$
for every $q, x > 0$ and $g : \r \to \r^+$ measurable, Theorem 2 (i) in \cite{R} yields
\begin{eqnarray}
\label{RE+}
\hU_q  f_M(x) & = & \Uapp_q  f_M(x)\; +\; \EE_x\!\lcr e^{-qT}\rcr\hU_q  f_M(0).
\end{eqnarray}
Recall from semigroup theory - see
e.g. Exercise VII (1.15) in \cite{RY} - that
$$\lim_{q\to \infty } q^2 \hU_q  f_M(x)\,-\,qf_M(x)\;= \; \hcR  f_M(x)\;= \; \hcR  f(x),$$
which altogether with the notation $\hT_1 =\inf\{t>0, \; \hZ_t \ge 1\}$ entails
\begin{eqnarray*}
  \hcR f(x) & = &  \lim_{q\to \infty} (q^2 \Uapp_q f_M(x)\,-\, q f_M(x)\; +\; q^2 \hU_q f_M(0) \EE_x[e^{-qT}])\\
& = &  \lim_{q\to \infty} (q^2 \Uapp_q f_M(x)\,-\, q f_M(x)\; +\; x^{-\aa} f(0) q \EE_0[e^{-q \hT_1}]),
\end{eqnarray*}
the second equality being a consequence of self-similarity and the easy fact $\lim_{q\rightarrow \infty } q \hU_q f_M(0)\;= \;f_M(0)\;= \;f(0).$
Proposition VIII.4 in \cite{B}, a standard self-similarity argument and the Tauberian theorem quoted in \cite{B} p. 10 give
$$\lim_{q\to \infty } q \,\EE_0[e^{-q \hT_1}]\; =\; \kappa$$
for some possibly vanishing, explicit constant $\kappa.$ This yields
$$\hcR f(x)\; =\;  \lim_{q\to \infty} (q^2 \Uapp_q f(x)\,-\, q f(x))\; +\; \frac{\kappa f(0)}{x^\aa}$$
and it remains to identify the limit on the right-hand side. Decomposing and changing the variable, one obtains
$$q^2 \Uapp_q f_M(x)\,-\, q f_M(x)\, =\, q^2\!\! \int_0^\infty \!\! e^{-qt} \EE_x [f_M(Z_{t})] dt \, - \, q f(x)\, -\, \int_0^\infty \!\! t e^{-t}\lpa \frac{\EE_x [f_M(Z_{t/q})\Un_{\{ T \le t/q\}}]}{t/q}\rpa dt.$$
By the same discussion as above, the  Markov property at time $T$ and the a.s. right-continuity of $t\mapsto Z_t$ at zero one has, recalling $f_M(x) = f(0)$ for every $-M <x < 0$,

$$\!\!\!\!\!\!\!\!\!\!\!\!\!\!\!\!\!\!\!\!\!\!\!\!\!\!\!\!\!\!\!\!\!\!\!\! \!\!\!\!\!\!\!\!\!\!\!\!\!\!\!\!\!\!\!\!\!\! \!\!\!\!\!\!\!\!\!\!\! \!\!\!\!\!\!\!\!\!\!\!    \frac{\kappa f(0)}{x^{\aa}} \;-\; \varepsilon(M)\;\le\;\liminf_{q\to +\infty} \frac{\EE_x [f_M(Z_{t/q})\Un_{\{ T \le t/q\}}]}{t/q}$$

\vspace{-2mm}

$$\;\;\;\;\;\;\;\;\;\;\;\;\;\;\; \;\;\;\;\;\;\;\;\;\;\;\;\;\;\;\;\;\;\;\;\;\;\;\;\;\;\;\;\;\;  \;\;\;\;\;\;\;\;\;\;\;\;\;\;\; \le\; \limsup_{q\to +\infty} \frac{\EE_x [f_M(Z_{t/q})\Un_{\{ T \le t/q\}}]}{t/q}\; \le\; \frac{\kappa f(0)}{x^{\aa}} \;+\; \varepsilon(M)$$
for every $t > 0$ and the {\em same} constant $\kappa$ as above, with $\varepsilon(M)\to 0$ as $M\to\infty.$ By Fatou's theorem, this entails

$$\!\!\!\!\!\!\!\!\!\!\!\!\!\!\!\!\!\!\!\!\!\!\!\!\!\!\!\!\!\!\!\!\!\!\!\! \!\!\!\!\!\!\!\!\!\!\!\!\!\!\!\!\!\!\!\!\!\! \!\!\!\!\!\!\!\!\!\!\! \!\!\!\!\!\!\!\!\!\!\!    \frac{\kappa f(0)}{x^{\aa}} \;-\; \varepsilon(M)\;\le\;\liminf_{q\to +\infty} \int_0^\infty t e^{-t}\lpa \frac{\EE_x [f_M(Z_{t/q})\Un_{\{ T \le t/q\}}]}{t/q}\rpa dt$$

$$\;\;\;\;\;\;\;\;\;\;\;\;\;\;\; \;\;\;\;\;\;\;\;\;\;\;\;\;\;\;\;\;\;\;\;\;\;\;\; \le\; \limsup_{q\to +\infty} \int_0^\infty t e^{-t}\lpa \frac{\EE_x [f_M(Z_{t/q})\Un_{\{ T \le t/q\}}]}{t/q}\rpa dt\; \le\; \frac{\kappa f(0)}{x^{\aa}} \;+\; \varepsilon(M).$$
One the other hand, again from the resolvent equation,
$$q^2 \int_0^\infty  e^{-qt} \EE_x [f_M(Z_{t})] dt \; - \; q f(x)\; =\;  q^2 \int_0^\infty  e^{-qt} \EE_x [f_M(Z_{t})] dt \; - \; q f_M(x)\; \to\; \cL f_M(x)$$
as $q\to +\infty,$ where $\cL$ is the infinitesimal generator of $Z.$ Indeed, one has $f_M\in\CC^2_b(\r)$ except possibly at zero where its left and right second derivatives are bounded and $f_M(x) \to 0$ as $\vert x\vert \to +\infty,$ so that $f_M\in \DD(\cL),$ as can be ssen readily from the proof of Theorem 31.5 in \cite{S}.
Supposing first $1 < \aa < 2,$ one has from the L\'evy-Khintchine formula
\begin{eqnarray*}
\cL f_M(x) & = & \int_\r (f_M(x+y) - f_M(x) - yf_M'(x))\nu(y) dy\\
& = & \int_0^x (f(x-y) - f(x) + yf'(x))\frac{c_-}{y^{\aa +1}} dy\; +\; \frac{c_-}{\aa x^{\aa}} (f(0) - f(x))\; +\; \frac{c_- f'(x)}{(\aa -1) x^{\aa-1}}\\
&  & \;\;\;\;\; \;\;\;\;\;\;\;\;\;\;\;\;\;\;\; +\;\; \int_0^\infty (f(x+y) - f(x) - yf'(x))\frac{c_+}{y^{\aa +1}} dy\; +\; \varepsilon (M).\\
\end{eqnarray*}

\vspace{-4mm}

\noindent
Letting $M\to +\infty,$ putting everything together and using the change of variable mentioned in \cite{BDP} involving the assumption $f\in\DD,$ one obtains
$$\hcR f(x)\; =\; \cL f(x)\; =\;\Gamma(-\aa) (c_-\Dma f (x) + c_+\Dpa f (x))\; +\; \frac{c_- f(0)}{\aa x^\aa}$$
as desired. The cases $\aa = 1$ and $0< \aa < 1$ are analogous and left to the reader.

\end{proof}

\begin{remarks} (a) The above constant $\kappa$ can be identified as
  $c_+/\aa,$ see Lemma 3.1 in \cite{CC} and the references therein. The value of this constant does not
  play any r\^ole here, but it is interesting to note that it is exactly the same as the one extracted from the L\'evy-Khintchine formula in the above proof.

\vspace{2mm}

(b) As in the proof of Proposition \ref{SMS}, it is possible to relax the condition $f\in\DD.$ For example when $\aa\in (1,2)$ the global boundedness condition on $f''$ can be changed into $(x^\gamma\wedge 1) f''(x)$ bounded on $(0,+\infty)$ for some $\gamma < 2-\aa,$ and when $\aa\in (0,1)$ the global boundedness condition on $f'$ can be changed into $(x^\gamma\wedge 1) f'(x)$ bounded on $(0,+\infty)$ for some $\gamma < 1-\aa.$ This is readily seen from the L\'evy-Khintchine formula and the proof of Theorem 31.5 in \cite{S}.

\vspace{2mm}

(c) With recurrent extensions, it is also possible to give an alternative proof to Proposition \ref{Dom}. Suppose as in Section 2 that $\aa\in (1,2)$ and that $Z$ has no positive jumps. By (\ref{RE+}), a function is in $\DD (\hL)$ iff it can be written
$$\hU_1  f(x)\; =\; \Uapp_1  f(x)\; +\; \EE_x [e^{-T}] \hU_1  f(0)$$
for some $f\in\CC_0.$  By self-similarity and Formula (7) in \cite{ThS1} - see also the references therein, we find first
$$\EE_x [ e^{-T}]\; =\;\EE_0 [ e^{- x^\aa \hT_1}]\; =\; \Fa(x) - \Fa'(x).$$
The term $\Uapp_1 f(x)$ can be handled with Suprun's formula. Specifically,
letting $a\to +\infty$ in Theorem 1 of \cite{B2} and using the discussion made
after Theorem 2 therein, we obtain
\begin{eqnarray*}
\Uapp_1 f(x) & = & \lpa\int_0^\infty e^{-y} f(y) dy\rpa\Fa'(x)\; -\;
 \int_0^x \Fa'(y) f(x-y) dy\\
& = & \lbd_f\Fa'(x)\; -\;
 \Fa'\star f \, (x)
\end{eqnarray*}
with the notation of Proposition \ref{Dom}. Last, we compute
$$\hU_1  f(0)\; =\; \EE\lcr \int_0^\infty e^{-t} f(Z_t - I_t) dt\rcr\; =\; \EE\lcr \int_0^\infty e^{-t} f(S_t) dt\rcr\; =\; \lbd_f$$
where the last equality follows from the discussion after (\ref{Deter}), paying here attention to the normalizing constants. Putting everything
together yields the expression for $\DD(\hL)$ given in Proposition
\ref{Dom}. The formula for $\DD(L)$ follows the same way, letting $a, x,
y\to +\infty$ with $a-x$ and $a-y$ constant in
Theorem 1 of \cite{B2} and identifying
$$\EE\lcr \int_0^\infty e^{-t} f(\hS_t) dt\rcr\; =\; \int_0^\infty (\Fa''(y)
-\Fa'(y)) f(y) dy.$$
We omit the details.
\end{remarks}

\subsection{Second proof of the theorem} In this paragraph we obtain a new proof of the identity (\ref{Cruz}) which does not depend on the results of \cite{ThS1} but on Mellin inversion. More precisely, we will show that
\begin{equation}
\label{Crux}
\EE_0[X_1^s]\; =\; \EE[\Va^s]\,\times\, \hEE_0[\hX_1^s]
\end{equation}
for every $s\in(1-\aa, \aa),$ which is plainly enough to get (\ref{Cruz}). We
start  with the fractional moments of the random variable $\Va,$ a computation
that could have been made directly by the residue theorem but since most of the argument was already carried out in \cite{ThS2} for some other purposes, we take the opportunity to shorten the proof.

\begin{lemma}
\label{plus}
For every $s\in(1-\aa,\aa),$ one has
$$\EE[\Va^s]\; =\;  \frac{\sin(\pi/\aa)\sin (\pi s)}{\aa\sin(\pi s/\aa)\sin(\pi(1-s)/\aa)}\cdot$$
\end{lemma}

\begin{proof} By equation (3) in \cite{ThS1} we know that the function
$$f_{\alpha}(t)\;=\;\frac{(- \sin\pi\aa)t^{\aa-1}(1+t)}{\pi(t^{2\alpha}-2t^{\aa}\cos\pi \aa +1)}$$
is a probability density over $\r^+.$ The fractional moments of the corresponding random variable $Y_\aa$ can be computed with the help of the beginning of the proof of Proposition 4 in \cite{ThS2} and a change of variable: one finds
$$\EE[Y_\aa^s]\; =\;  \frac{\sin(\pi/\aa)\sin (\pi s)}{\aa\sin(\pi s/\aa)\sin(\pi(1+s)/\aa)}$$
for every $s\in(-\aa, \aa-1).$ Notice that making $s = -1$ entails
$$\int_0^\infty v_{\aa}(t) \, dt\; =\; \EE[Y_\aa^{-1}]\; =\; 1,$$
which shows that $v_\aa$ is a probability density with an argument slightly different from the introduction. Finally, the fractional moments of $\Va$ are given by
$$\EE[\Va^s]\; =\; \EE[Y_\aa^{s-1}]\; =\;  \frac{\sin(\pi/\aa)\sin (\pi s)}{\aa\sin(\pi s/\aa)\sin(\pi(1-s)/\aa)}$$
for every $s\in(1-\aa,\aa).$

\end{proof}

To compute the fractional moments of $X_1$ under $\pb_0,$ we will need more material on recurrent extensions and exponential functional of L\'evy processes. With the notations of Section 2, let $\hZ = -Z$ be the dual process and $\hpb_x$ its law starting at $x >0.$ Introducing the stopping time $\hT = \inf \{s>0,\: \hZ_s<0\},$ consider the positive $(1/\aa)-$self-similar Feller process
$$\hR_t \; =\; \hZ_t \Un_{\{ \hT > t\}}.$$
The well-known Lamperti transformation \cite{L} shows that the process defined $\xi_t = \log \hR_{\tau_t}$ for every $t< \hT,$ with the notation
$$\tau_t\;=\;\inf\lacc s>0,\; \int_0^s \hR_r^{-\alpha}dr>t\racc,$$
is a L\'evy process starting at $\log x $. Its L\'evy-Khintchine exponent $\psi$ which is defined by
$$\EE[e^{-\lbd \hxi_t}] \; =\; e^{t \psi (\lbd)}$$
for every $\lbd \ge 0$ - recall that $\xi$ has no negative jumps so that the above expectation is finite, has been computed in  \cite{CC} in terms of a certain improper integral. The next lemma gives a more tractable formulation in terms of Gamma functions.

\begin{lemma}
\label{LK}
With the normalization of Section 2, one has $\psi(\lbd) =\Gamma(\lbd+\alpha)/\Gamma(\lbd).$
\end{lemma}

\begin{proof} By Corollary 1 in \cite{CC} and Theorem 2.4 of \cite{Pa1}, one has
\begin{eqnarray*}
\psi(\lbd)\; =\; \Phi (\ii\lbd) & = & \frac{\lbd}{(\aa -1)\Gamma(-\aa)}\; +\; \int^{+\infty}_0 (e^{-\lbd y} - 1 + \lbd (e^y -1)\Un_{\{e^y \le 2\}})\frac{e^y dy}{\Gamma(-\aa) (e^y - 1)^{1+\aa}}\\
& = & \int_1^{+\infty} \lpa \frac{y^{-\lbd} - 1 + \lbd (y -1)}{\Gamma(-\aa) (y - 1)^{1+\aa}}\rpa dy\\
& = &  \int_0^1 \lpa \frac{(y^{\lbd} - 1) y^{\aa -1} - \lbd (y -1)}{\Gamma(-\aa) (1-y)^{1+\aa}}\rpa dy\; +\; \lbd \int_0^1 \lpa \frac{y^{\aa -2} - 1}{\Gamma(-\aa) (1-y)^\aa}\rpa dy\\
& = & \frac{\Gamma(\lbd+\alpha)}{\Gamma(\lbd)}\; -\; \frac{\lbd}{(\aa -1)\Gamma(-\aa)}\; +\; \lbd \int_0^1 \lpa \frac{y^{\aa -2} - 1}{\Gamma(-\aa) (1-y)^\aa}\rpa dy\\
\end{eqnarray*}
where we have used several changes of variable. The last integral can be computed with the help of Formula (2.3) in \cite{Pa1}: one gets
$$\int_0^1 \lpa \frac{y^{\aa -2} - 1}{(1-y)^\aa}\rpa dy\; =\; -\Gamma(1-\aa) (2-\aa)_\aa \; =\; 1/(\aa-1),$$
yielding the desired formula for $\psi.$
\end{proof}

\begin{remark} Supposing now that $Z$ has positive jumps with the notations of
  Paragraph 3.1, an analogous simplification of Corollary 1 in \cite{CC} with
  Theorem 2.4 of \cite{Pa1} shows the following general formula for the
  L\'evy-Khintchine exponent of the L\'evy process associated to $\hR$ by the
  Lamperti transformation:
$$\psi(\lbd) \; =\; \Gamma(-\aa)(c_-\Gamma(\lbd+\alpha)/\Gamma(\lbd) + c_+
\Gamma(1-\lbd)/\Gamma(1- (\aa +\lbd)))$$
for every $\lbd\in (-\aa, 1).$ Setting $\theta = \inf\{\lbd > 0: \;
\psi(-\lbd) =0\},$ a simple analysis shows
then that $\theta =\aa\rho = \aa(1-{\hat \rho})$ where ${\hat \rho}$ is the
asymmetry coefficient of $\hZ,$ which can also be checked in considering the
invariant function of $\hR$ - see again Example 3 in \cite{R}.

\end{remark}

Setting now
$$\thp \; =\; \inf\{\lbd > 0: \; \psi(-\lbd) =0\}\; =\;1\; <\; \aa,$$
Theorem 2 in \cite{R} shows the existence of a unique recurrent extension for $\hR$ leaving 0 continuously, whose resolvent $(\Uam_q)_{q\ge 0}$ is characterized by the formula
\begin{eqnarray}
\label{REP}
\Uam_q f(0)& = & \frac{\aa q^{-1/\aa}}{\Gamma(1-1/\aa)\EE [ I_-^{1/\aa-1}]}\int_0^{\infty}f(y) y^{\aa-2}\EE[e^{-qy^\aa I_-}] dy
\end{eqnarray}
for any positive measurable function $f,$ where
\begin{equation}
\label{IP}
I_- = \int_0^{\infty} e^{\xi^-_s}ds
\end{equation}
and $\xi^-$ is a spectrally negative L\'evy process with Laplace exponent
$$\psip(\lbd)\; =\;\psi(\aa\lbd-\thp)\; =\; \;\psi(\aa\lbd-1)\; =\;\frac{\Gamma(\aa(\lbd +1) -1)}{\Gamma(\aa\lbd -1)}\cdot$$
As mentioned in Example 3 of \cite{R}, this recurrent extension is the process $X,$ and with the help of this identification we can now compute the fractional moments of $X_1$:

\begin{prop}
\label{minus}
For any $s\in (1-\aa, \aa)$, one has
$$\EE_0[X_1^s]\;=\; \frac{\sin (\pi/\aa)\sin (\pi s) \Gamma (s+1)}{\alpha \Gamma (s/\aa +1)
\sin (\pi s/\aa) \sin (\pi (1-s)/\aa)}\cdot$$
\end{prop}
\begin{proof} Let us first connect the moments of $X_1$ under $\pb_0$ to those of $I_-.$ Introducing the positive measurable function $p_s (t) = t^s$ over $\r^+,$ for every $s\in
(1-\aa, \aa)$ one has
\begin{eqnarray*}
\EE_0[X_1^s] & = & \frac{q^{s/\aa+1}}{\Gamma(s/\aa+1)}
\int_0^\infty e^{-qt}t^{s/\alpha}\EE_0[X_1^s] dt\\
 & = & \frac{q^{s/\aa+1}}{\Gamma(s/\aa+1)} \int_0^\infty e^{-qt} \EE_0[X_t^s] dt \; =\; \frac{q^{s/\aa+1} U_q^-p_s (0)}{\Gamma(s/\aa+1)}
 \end{eqnarray*}
where the second equality comes from self-similarity and the third from the fact that the resolvent of $X$ is $(\Uam_q)_{q\ge 0}$. From (\ref{REP}) and after some simplifications, this entails
$$\EE_0[X_1^s] \; =\; \frac{\Gamma (1+(s-1)/\aa)\EE[I_-^{-(1 +(s-1)/\aa)}]}{\Gamma(1-1/\aa)\Gamma(s/\aa+1)\EE[I_-^{1/\aa-1}]}$$
for any $s\in (1-\aa, \aa).$  On the other hand, from Theorem 2.1 and Formula (2.1) in \cite{Pa3}  - with our notation which entails  $\gamma = 1/\aa$ therein, see \cite{Pa4} for details - we know that the density function $f_-$ of $I_-$ has the alternate series representation
$$f_-(t)\;= \;C_\aa \sum_{n=0}^{\infty} (-1)^n \lpa\frac{\Gamma (n+1 +1/\aa)}{\Gamma (\aa(n+1))} \rpa t^{-n+1+1/\aa}, \qquad t>0,$$
where $C_\aa$ is a positive constant to be determined below. This
representation of the density prevents from computing the fractional moments
of $I_-$ by direct integration. Instead, one can use a so-called Mellin-Barnes
integral representation of $f_-$, which is obtained simply after a contour
integration along a big half-circle in the half-plane $x > -1$ - see e.g. Section 3.4 in \cite{KP} for details. For any $c\in (-1,0)$ one has
 \begin{eqnarray*}
 f_-(t) &=& \frac{C_\aa}{2\pi \ii} \int_{c-\ii\infty}^{c+\ii\infty} \lpa\frac{\Gamma (s+1+1/\aa)\Gamma(s+1)\Gamma(-s)}{\Gamma(\aa(s+1))}\rpa t^{-s+1+1/\aa}\, ds  \\
 &=& \frac{C_\aa}{2\pi \ii}  \int_{c'-\ii\infty}^{c' +\ii\infty} \lpa\frac{\Gamma(s)\Gamma(s- 1/\aa)\Gamma(-s+1 +1/\aa)}{\Gamma\left(\alpha s-1\right)}\rpa t^{-s} ds
\end{eqnarray*}
after a change of variable and taking any $ c'\in (1/\aa, 1+1/\aa).$ The inversion formula for the Mellin transform -  see e.g. formula (3.1.5) in \cite{KP} - entails then
\begin{eqnarray*}
\EE[I_-^{s}] & = & C_\aa\frac{\Gamma (s+1)\Gamma(s+1-1/\aa)\Gamma(-s+ 1/\alpha)}{\Gamma (\aa (s+1)-1)}
\end{eqnarray*}
for every $s\in (-1+1/\aa, 1/\aa).$ Notice in passing, though we shall not
need this, that making $s=0$ allows also to compute $C_\aa = \Gamma(\aa-1)/(\Gamma(1-1/\aa)\Gamma(1/\aa)).$ After some simple transformations, we finally deduce that for any $s\in (1-\aa, \aa)$,
$$ \EE_0[X_1^s] \;= \; \frac{\Gamma(1+ (s-1)/\aa)\Gamma((1-s)/\aa)\Gamma(-s/\aa)}{\aa\Gamma (1-1/\aa)\Gamma(1/\aa)\Gamma(-s)}\; =\;\frac{\sin (\pi/\aa)\sin (\pi s) \Gamma (s+1)}{\alpha \Gamma (s/\aa +1) \sin (\pi s/\aa) \sin (\pi (1-s)/\aa)}\cdot$$
\end{proof}

\noindent
{\bf End of the proof.} The property that the law of $\hX_1$ under $\hpb_0$ is that of $S_1\elaw T_1^{-1/\aa}$ under $\pb$ and a well-known, aforementioned moment formula for $T_1$ entails
$$\hEE_0[\hX_1^s]\;=\; \frac{\Gamma (s+1)}{\Gamma (s/\aa +1)}$$
for every $s > -1,$ so that (\ref{Crux}) simply follows from Lemma \ref{plus} and Proposition \ref{minus}.

\vspace{2mm}

\noindent
{\bf Acknowledgements.} The research of PP was supported by a grant from the
National Bank of Belgium. TS wishes to thank Grant ANR-09-BLAN-0084-01.

\end{document}